\documentclass{amsart}
\usepackage{amssymb,amsmath,graphicx,mathrsfs,amsfonts}
\usepackage{tikz}
\usepackage{url}
\usepackage[dvips]{epsfig}
\newtheorem{thrm}{Theorem}[section]
\newtheorem{lem}[thrm]{Lemma}
\newtheorem{prop}[thrm]{Proposition}

\newtheorem{remark}[thrm]{Remark}
\usepackage{graphicx}
\DeclareUnicodeCharacter{0301}{\'{e}}
\begin{document}
\author[L.~G.~Molinari]{Luca Guido Molinari}
\address{L. G. Molinari, Physics Department ``Aldo Pontremoli'',
Universit\`a degli Studi di Milano and I.N.F.N. sez. di Milano,
Via Celoria 16, 20133 Milano, Italy.}
\email{luca.molinari@unimi.it}
\begin{abstract}
I conjecture three identities for the determinant of adjacency matrices of graphene triangles and trapezia 
with Bloch (and more general) boundary conditions. 
For triangles, the parametric determinant is equal to the characteristic polynomial of the symmetric Pascal matrix. 
For trapezia it is equal to the determinant of a sub-matrix. 
Finally, the determinant of the tight binding matrix equals its permanent. 
The conjectures are supported by analytic evaluations and Mathematica, for moderate sizes.
They establish connections with counting problems of partitions, lozenge tilings of hexagons, dense loops on a cylinder. 
\end{abstract}
\title[Graphene nanocones and Pascal matrices] 
{GRAPHENE NANOCONES  and PASCAL MATRICES\\
\Small{\it -- some determinantal conjectures --}}
\maketitle


\section{Introduction}
Physics often entails intriguing mathematics, and Carbon Nanocones are no exception.
They are conical honeycomb lattices of Carbon atoms, with a cap of Carbon rings causing different cone apertures \cite{Balaban94}. 
The ones with smallest angle, near $19^\circ$, were sinthetized by Ge and Sattler in 1994, in the hot vapour phase of carbon; the others were found by accident, in pyrolysis of heavy oil \cite{Sattler08}. 
\begin{figure}[h] 
   \centering
   \includegraphics[width=2.5in]{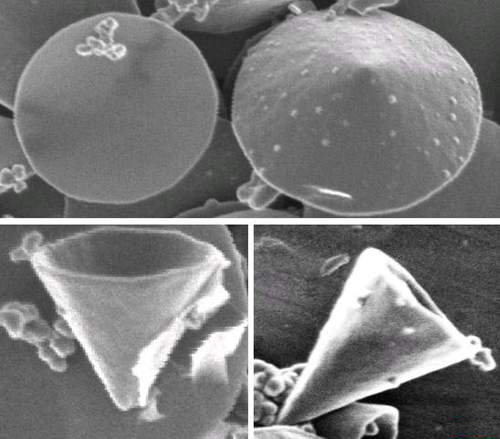} 
   \caption{{\footnotesize Scanning electron microscope images of carbon nanocones (size order $1\mu m$) produced 
   at the Kvaerner Carbon Black \& Hydrogen Process \cite{Knudsen}.}}
   \label{fig:1}
\end{figure}

The simplest description for such structures is through the adjacency matrix. Chemists 
give weight $-\beta $ to lattice links among atoms, and name it H\"uckel matrix. 
Because of the $C_n$ symmetry of nanocones, the H\"uckel matrix can be restricted to a single 
honeycomb triangle, with Bloch boundary conditions $e^{\pm i 2\pi/n}$ on two sides. 
Properties of the spectrum, based on graph theory and symmetries, were studied in \cite{Heiberg07}.
Recently, based on numerics, new rules were established for the filling of levels, that differ from the 
ordinary H\"uckel rules \cite{Evan1}. \\
This paper originates from a question by an author in \cite{Evan1}: what is the determinant of the Huckel
matrix of a nanocone? The answer here conjectured opens a nice connection with Pascal matrices and combinatorics.

\begin{figure}[h]
\begin{center}
\includegraphics[width=4.5cm,clip=]{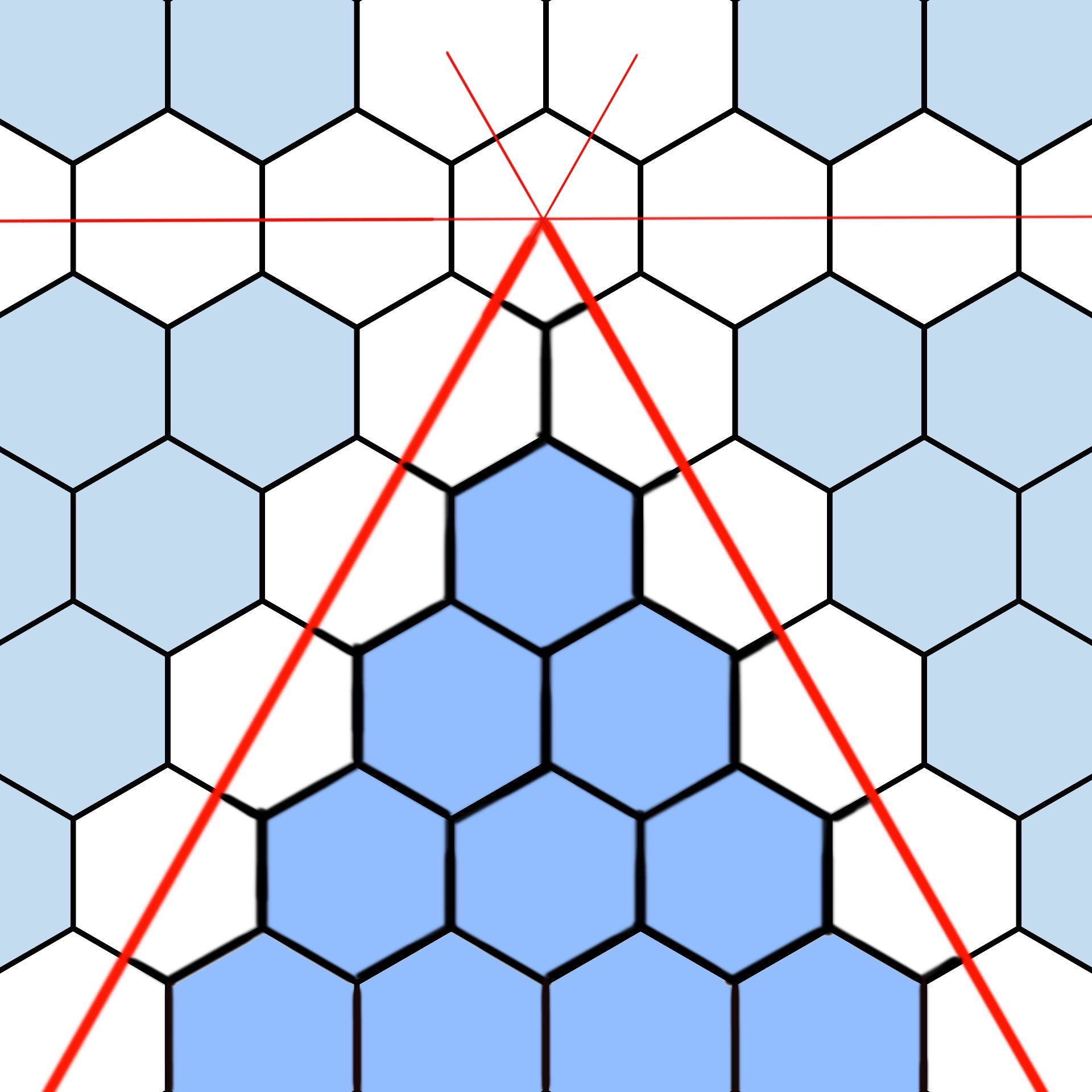}
\caption{\footnotesize Six lines cut graphene into triangles 
with zig-zag rows of 1, 3, 5, ... atoms. Cones are obtained by re-joining
3, 4, 5 triangles, with a cap-ring of 3, 4, 5 atoms.} 
\end{center}
\end{figure}
It is convenient to replace the Bloch factors by free parameters $x$ and $y$, or even by parameters $x_i$ and $y_i$ for each row.\\
The size of a H\"uckel matrix with $n$ rows is $(n+1)^2$, which is the number of atoms in the triangle. \\
\begin{minipage}{.4 \textwidth}
\begin{tikzpicture}
\draw[densely dotted,thick] (-2,-1) -- (-2.5,-1);
\draw[densely dotted,thick] (2,-1) -- (2.5,-1);
\draw[densely dotted,thick] (-1,0) -- (-1.5,0);
\draw[densely dotted,thick] (1,0) -- (1.5,0);
\draw[densely dotted,thick] (-0.5,1) -- (0.5,1);
\draw[black, thick] (-2,-1) -- (2,-1);
\draw[black, thick] (-1,-1) -- (1,-1);
\draw[gray, thick] (-1,-1) -- (-1,0);
\draw[gray, thick] (1,-1) -- (1,0);
\draw[black, thick] (-1,0) -- (1,0);
\draw[gray, thick] (0,0) -- (0,1);
\filldraw[blue] (0,1) circle (2pt) node[anchor=south] {0};
\filldraw[blue] (1,0) circle (2pt) node[anchor=north] {};
\filldraw[red] (0,0) circle (2pt) node[anchor=north] {};
\filldraw[blue] (-1,0) circle (2pt) node[anchor=south] {1};
\filldraw[blue] (2,-1) circle (2pt) node[anchor=north] {};
\filldraw[red] (1,-1) circle (2pt) node[anchor=north] {};
\filldraw[blue] (0,-1) circle (2pt) node[anchor=north] {};
\filldraw[red] (-1,-1) circle (2pt) node[anchor=north] {};
\filldraw[blue] (-2,-1) circle (2pt) node[anchor=south] {2};
\end{tikzpicture}
\end{minipage}
\begin{minipage}{.1 \textwidth}
\begin{align*}
\small{H_2=
\left [ \begin{array}{c|ccc|ccccc}
S_0 & 0 &1 & & & & & & \\
\hline
0& 0& 1& y_1& 0& 1& & & \\
1& 1& 0& 1& 0& & 0& & \\
& x_1& 1& 0& 0& & & 1& \\
\hline
& 0&0 &0 & 0& 1& & & y_2\\
& 1& & & 1& 0& 1& & \\
& & 0& & & 1& 0& 1& \\
& & & 1& & & 1& 0& 1\\
& & & & x_2& & & 1& 0
\end{array}\right] }
\end{align*}
\end{minipage}
{}\\\\\\
{\footnotesize Example. A triangle with one hexagon and its H\"uckel matrix.
In this representation, the diagonal blocks 
are the rows of the graph, with 1,3,5 atoms. Edges are 1s in the matrix. Dots connect row extrema with weights $x_i=y_i=1$  ($S_0=x_0+y_0$).}\\ 

The determinants of small H\"uckel matrices are evaluated:
\begin{align*}
\det H_0=&x+y\\
\det H_1=&x^2+3xy+y^2\\ 
\det H_2=&x^3 + 9 x^2 y + 9 x y^2 + y^3\\
\det H_3=&x^4 + 29 x^3 y + 72 x^2 y^2 + 29 x y^3 + y^4\\
\det H_4 =&x^5 + 99 x^4 y + 626 x^3 y^2 + 626 x^2 y^3 + 99 x y^4 + y^5\\
\det H_5 =&x^6 + 351 x^5 y + 6084 x^4 y^2 + 13869 x^3 y^3 + 6084 x^2 y^4 +  351 x y^5 + y^6\\
\det H_6=& x^7 + 1275 x^6 y + 64974 x^5 y^2 + 347020 x^4 y^3 + 347020 x^3 y^4  \\
&\quad + 64974 x^2 y^5 + 1275 x y^6 + y^7
 \end{align*}
$\bullet $ The determinants are homogeneous palindromic polynomials with positive coefficients.\\
$\bullet $ If the parameters are $x_1,x_2,... ,y_n$ are kept different, the determinants are homogeneous polynomial of degree $n+1$ in the variables, with coefficients that are perfect squares of integers (this feature was noticed to me by Andrea Sportiello). \\
 
Their evaluation for small size $n$ with Schur's formula for block matrices showed that the determinant of a H\"uckel matrix of size $(n+1)^2$ eventually reduces to the determinant of a matrix of binomials of size $(n+1)$ (a step is illustrated in section 8). Binomials occur in many counting problems of benzenoids \cite{Bodroza88, Doslic07}.\\
After discovering the beautiful Pascal matrices in Lunnon's paper \cite{Lunnon77}, I found with Mathematica that the polynomial coefficients in $\det H_n$ for accessible values of $n$, match the sequences A045912 in OEIS \cite{OEIS} for the characteristic polynomials of Pascal matrices. This led me to the first conjecture; the others then came out. 

Conjectures 1 and 2 imply symmetries of the eigenvalues of H\"uckel matrices of nanocones and the larger family of graphannulenes, that support the generalized H\"uckel rules for the electronic structure proposed by Evangelisti et al. in ref.\cite{Evan1}.

Before stating the conjectures, let's introduce Pascal matrices.

 \section{Pascal matrices}
There are various forms of Pascal matrices, and much literature \cite{Lunnon77,Brawer92,Yates14,Edelman04}. The lower triangular Pascal matrix $P_n$ has size $n+1$ and elements that are the binomial coefficients of the Pascal triangle. The inverse matrix also has 
binomial coefficients, but with alternating signs. For example: 
 \begin{align*}
 P_4 =
 \left[ \begin{array}{ccccc}
 1 &  & & & \\
1 &1 &  & &\\
 1 & 2 & 1& &\\
1 & 3& 3&1&\\
1& 4 &6 & 4  & 1 
 \end{array}\right ] , \quad
 P_4^{-1} = \left[ \begin{array}{ccccc}
 1 &  & & & \\
-1 &1 &  & &\\
 1 & -2 & 1& &\\
-1 & 3& -3&1&\\
1& -4 &6 & -4  & 1 
 \end{array}\right ]
 \end{align*} 
 The product $Q_n=P_nP_n^T$ is the symmetric positive Pascal matrix. It has unit determinant and matrix elements $ (Q_n)_{ij} = \binom{i+j}{j}$ $i,j=0,...,n $:
 
 \begin{align*}
 Q_4 = \left[ \begin{array}{ccccc}
 1 & 1 &1 &1 &1 \\
1 &2 & 3 &4 & 5\\
 1 & 3 & 6& 10& 15\\
1 & 4& 10 & 20 & 35\\
1& 5 &15 & 35  & 70 
 \end{array}\right ]
 \end{align*}
 Since the eigenvalues are strictly positive, 
 the coefficients of the polynomial $\det (z+Q_n)$ are positive. Moreover the eigenvalues come in pairs: $q_{n+1-j} = 1/q_j$;  
 if the matrix size is odd ($n$ even), an eigenvalue of $Q_n$ is 1.  \\
 Note that the entry $(i,j)$ (independent of $n$) counts the number of paths with $i+j$ steps $\rightarrow$ or $\downarrow$ that connect 
 the corner ($0,0$) with ($i,j$). For example, ($3,4$) can be reached in $35=\binom{7}{3}$ different ways from ($0,0$).

\section{The first conjecture (honeycomb triangles)}
A graphene triangle and its H\"uckel matrix are shown below. For later convenience the atoms are marked as red and blue (two triangular 
sub-lattices of hexagonal graphene). \\
\begin{minipage}{.4 \textwidth}
\small{
\begin{tikzpicture}
\draw[densely dotted,thick] (-2,2.5) -- (-2.5,2.5);
\filldraw[blue] (-2, 2.5) circle (2pt) node[anchor=south] {4};
\draw[black, thick] (-2,2.5) -- (-1.5,3);
\filldraw[red] (-1.5,3) circle (2pt) node[anchor=north] {};
\draw[black, thick] (-1.5,3) -- (-1,2.5);
\filldraw[blue] (-1,2.5) circle (2pt) node[anchor=north] {};
\draw[black, thick] (-1,2.5) -- (-0.5,3);
\filldraw[red] (-0.5,3) circle (2pt) node[anchor=north] {};
\draw[black, thick] (-0.5,3) -- (0,2.5);
\filldraw[blue] (0,2.5) circle (2pt) node[anchor=north] {};
\draw[black, thick] (0,2.5) -- (0.5,3);
\filldraw[red] (0.5,3) circle (2pt) node[anchor=north] {};
\draw[black, thick] (0.5,3) -- (1,2.5);
\filldraw[blue] (1,2.5) circle (2pt) node[anchor=north] {};
\draw[black, thick] (1,2.5) -- (1.5,3);
\filldraw[red] (1.5,3) circle (2pt) node[anchor=north] {};
\draw[black, thick] (1.5,3) -- (2,2.5);
\filldraw[blue] (2,2.5) circle (2pt) node[anchor=north] {};
\draw[densely dotted,thick] (2,2.5) -- (2.5,2.5);
%
\draw[gray, thick] (-1.5,3) -- (-1.5,3.5);
\draw[gray, thick] (-0.5,3) -- (-0.5,3.5);
\draw[gray, thick] (0.5,3) -- (0.5,3.5);
\draw[gray, thick] (1.5,3) -- (1.5,3.5);
\draw[densely dotted,thick] (-1.5,3.5) -- (-2,3.5);
\filldraw[blue] (-1.5,3.5) circle (2pt) node[anchor=south] {3};
\draw[black, thick] (-1.5,3.5) -- (-1,4);
\filldraw[red] (-1,4) circle (2pt) node[anchor=north] {};
\draw[black, thick] (-1,4) -- (-0.5,3.5);
\filldraw[blue] (-0.5,3.5) circle (2pt) node[anchor=north] {};
\draw[black, thick] (-0.5,3.5) -- (0,4);
\filldraw[red] (0,4) circle (2pt) node[anchor=north] {};
\draw[black, thick] (0,4) -- (0.5,3.5);
\filldraw[blue] (0.5,3.5) circle (2pt) node[anchor=north] {};
\draw[black, thick] (0.5,3.5) -- (1,4);
\filldraw[red] (1,4) circle (2pt) node[anchor=north] {};
\draw[black, thick] (1,4) -- (1.5,3.5);
\filldraw[blue] (1.5,3.5) circle (2pt) node[anchor=north] {};
\draw[densely dotted,thick] (1.5,3.5) -- (2,3.5);
\draw[gray, thick] (-1,4) -- (-1,4.5);
\draw[gray, thick] (0,4) -- (0,4.5);
\draw[gray, thick] (1,4) -- (1,4.5);
\draw[densely dotted,thick] (-1,4.5) -- (-1.5,4.5);
\filldraw[blue] (-1,4.5) circle (2pt) node[anchor=south] {2};
\draw[black, thick] (-1,4.5) -- (-0.5,5);
\filldraw[red] (-0.5,5) circle (2pt) node[anchor=north] {};
\draw[black, thick] (-0.5,5) -- (0,4.5);
\filldraw[blue] (0,4.5) circle (2pt) node[anchor=south] {};
\draw[black, thick] (0,4.5) -- (0.5,5);
\filldraw[red] (0.5,5) circle (2pt) node[anchor=north] {};
\draw[black, thick] (0.5,5) -- (1,4.5);
\filldraw[blue] (1,4.5) circle (2pt) node[anchor=south] {};
\draw[densely dotted,thick] (1,4.5) -- (1.5,4.5);
\draw[gray, thick] (-0.5,5) -- (-0.5,5.5);
\draw[gray, thick] (0.5,5) -- (0.5,5.5);
\draw[densely dotted,thick] (-0.5,5.5) -- (-1,5.5);
\filldraw[blue] (-0.5,5.5) circle (2pt) node[anchor=south] {1};
\draw[black, thick] (-0.5,5.5) -- (0,6);
\draw[black, thick] (0.5,5.5) -- (0,6);
\filldraw[blue] (0.5,5.5) circle (2pt) node[anchor=north] {};
\draw[densely dotted,thick] (0.5,5.5) -- (1,5.5);
\filldraw[red] (0,6) circle (2pt) node[anchor=north] {};
\draw[gray, thick] (0,6) -- (0,6.5);
\filldraw[blue] (0,6.5) circle (2pt) node[anchor=south] {0};
\draw[densely dotted,thick] (-0.5,6.5) -- (0.5,6.5);
\end{tikzpicture}
}
\end{minipage}
\begin{minipage}{.\textwidth}
{\small
\begin{align*}
H_5({\bf x},{\bf y})=\left [\begin{array}{cccc}
 T_0& R_1^T &  &\\
 R_1& \ddots& \ddots &\\
 &\ddots&\ddots &R_4^T\\
&& R_4& T_4
\end{array}\right] 
\end{align*}
}
\end{minipage}\\
A block $T_m(x_m,y_m)$ is a square matrix of size $2m+1$ that describes a row of $2m+1$ atoms with
boundary parameters $x_m,y_m$: $T_0(x_0,y_0)=x_0+y_0$, 
\begin{align*}
{\small, \quad 
T_1 (x_1,y_1) = \left [ \begin{array}{ccc} 0 & 1 & y_1 \\ 1 & 0 & 1  \\ x_1 & 1 & 0  \end{array}\right ], 
\quad
T_2 (x_2,y_2) = \left [ \begin{array}{ccccc} 0 & 1 & 0 & 0 & y_2 \\ 1 & 0 & 1 & 0 & 0\\ 0 & 1 & 0 & 1& 0\\ 0 & 0 &1& 0 & 1 \\ 
x_2 & 0 & 0 & 1 & 0 \end{array}\right ], \ldots 
}
\end{align*}
\noindent
The size of $H_n$ is $1+3+...+(2n+1)=(n+1)^2$. 
The blocks $R_m$ are $(2m+1)\times (2m-1)$, with unit elements for the vertical edges in the graph, 
joining atoms in rows $m-1$ and $m$:
\begin{align*}{\small
R_1 =\left [\begin{array}{c}
  0 \\
  1 \\
  0 
\end{array}\right] , \quad 
R_2 =\left [\begin{array}{ccc}
  0& 0 & 0 \\
  1& 0&0 \\
 0 & 0&0 \\
0  &0 & 1 \\
  0& 0& 0
\end{array}\right] , \quad 
R_3 =\left [\begin{array}{ccccc}
0&0& 0&0&0 \\
1& 0&0 &0&0\\
0 & 0&0 &0&0\\
0  &0 & 1 &0&0\\
0&0&  0& 0& 0\\
0  &0 & 0 &0&1\\
0  &0 & 0 &0&0
\end{array}\right], \quad \ldots }
\end{align*}

\noindent
{\bf Conjecture 1:}\label{CJ1}
{\em The determinant of $H_n({\bf x},{\bf y})$ is equal to the determinant of the following matrix of size $n+1$:}
\begin{gather}\label{333}
 \left [ \begin{array}{cccccc}
x_n+y_n & -\binom{n}{1} y_n & \binom{n}{2} y_n & -\binom{n}{3}y_n & ... & \pm \binom{n}{n}y_n \\
\binom{n}{1}x_n & x_{n-1}+y_{n-1} & -\binom{n-1}{1}y_{n-1}  & \binom{n-1}{2} y_{n-1} & ... & \mp \binom{n-1}{n-1}y_{n-1}\\
\binom{n}{2}x_n &\binom{n-1}{1}x_{n-1} & x_{n-2}+y_{n-2} & - \binom{n-2}{1}y_{n-2}& ... & \vdots\\
\vdots &\vdots & \vdots & ... & ... & \vdots \\
\vdots &\vdots &\vdots & ... & x_1 +y_1& -\binom{1}{1}y_1 \\
\binom{n}{n} x_n &\binom{n-1}{n-1}x_{n-1} & \binom{n-2}{n-2} x_{n-2} & ... & \binom{1}{1}x_1& x_0+y_0\\
\end{array} \right ]. \nonumber
\end{gather}
\\
Example: 
\begin{align*}
&\det H_3= \left | \begin{array}{cccc}
 x_3+y_3 & -3y_3& 3y_3 & -y_3\\
 3x_3& x_2+y_2 & -2y_2  & y_2\\
 3x_3 & 2x_2 & x_1 +y_1& -y_1 \\
 x_3 & x_2& x_1& x_0+y_0
\end{array} \right |\\
&=(x_0+y_0) [x_1 x_2 x_3 + x_2 x_3 y_1 + x_1 x_2 y_3+ x_1 y_2 y_3+ x_1 x_3 y_2 + x_3 y_1 y_2 + x_2 y_1 y_3
 +y_1 y_2 y_3\\
&+4 x_2 y_2 (x_3+y_3)+9 x_3 y_3 (x_1+y_1+ x_2 +y_2)]
+ x_1 x_2 x_3 (y_1+ y_2+y_3) +  (x_1 +x_2)x_3 y_1 y_2 \\ 
& + x_1 x_2 (y_1+y_2) y_3    + x_2 x_3 y_2 y_3  
+ (x_1+x_2 +x_3)y_1 y_2 y_3  + 4 x_3 y_3( x_1 y_2 + x_2 y_1)+9 x_1 x_3 y_1 y_3
\end{align*}
It is a sum of monomials of degree $4$ with coefficients that are perfect integer squares.\\

With all $x_k=x$ and $y_k=y$, we obtain the interesting identity:
\begin{prop}
\begin{align}
\boxed{ \det H_n(x,y) = x^{n+1} \det \left (Q_n +\frac{ y}{x} I_{n+1} \right ) } \label{char}
\end{align}
where $Q_n$ is the symmetric Pascal matrix.
\begin{proof}
The $(n+1)\times (n+1)$ matrix in \eqref{333} is the sum of  lower and upper triangular matrices related to the Pascal matrix: 
$ x (J_nP_n^TJ_n)  + y (J_nP_n^{-1} J_n) $, where 
$$J_n= \left(\begin{smallmatrix} & &  1\\ & \reflectbox{$\ddots $}  &\\ 1 &  &  \end{smallmatrix} \right) $$ 
is the inversion matrix. Then: 
\begin{align*}
 \det H_n(x,y) &= \det (x P_n^T + y P_n^{-1}) \\ &= \det (xP_n ^{-1}) \det ( P_n P_n^T+y/x)\\
&=x^{n+1}\det(Q_n+y/x)
\end{align*}
More generally: $\det H_n({\bf x},{\bf y}) = \det (J_nP_n^TJ_nX +  YJ_nP_n^{-1}J_n) = \det(P_n J_nXJ_n +J_nYJ_n P_n^{-1})$ where $X={\rm diag}[x_i]$ and $Y={\rm diag}[y_i]$. If $Y$ is the identity matrix,
then $\det H_n({\bf x},1) = (x_0...x_n) \det (Q_n + {\rm diag}(1/x_0,...,1/x_n))$.
For example: the size of $H_3$ is
$N=16$, and
$$\det H_3 ({\bf x},1)  = (x_0...x_3) \det \left[ \begin{array}{cccc} 
1+1/x_0 & 1& 1& 1\\
 1& 2+1/x_1 & 3  & 4\\
 1& 3 & 6 +1/x_2& 10 \\
 1& 4& 10& 20+1/x_3 \end{array} \right ]. $$
\end{proof}
\end{prop}
In the next section I report the very few known analytic values of $\det H_n$, given by eq.\eqref{char} for special values of
the ratio $\omega=y/x=e^{2i\theta}$. They result from amazing combinatorial problems, that are here reported for the delight of the reader.

\section{Pascal matrix, partitions and lozenge tilings}
The characteristic polynomial of the symmetric Pascal matrix appears in diverse combinatorial problems 
with the following generalization:
\begin{align}
\det \left [\binom{m+j+k}{k} + \omega \delta_{jk}  \right ]_{0\le j,k\le n} \quad m=0,1,2, \ldots  \label{dets}
\end{align}
George Andrews \cite{Andrews} in 1979, in proving the weak Macdonald conjecture for certain plane partitions, 
obtained the closed expression of \eqref{dets} with $\omega=1$. He exploited identities for
hypergeometric series.
For the Pascal matrix ($m =0$) it greatly simplifies:
\begin{align}
\det (Q_n+I_{n+1})  =  \prod_{k=0}^n \frac{1}{3k+1} \frac{k! (3k+2)! }{(2k)! (2k+1)!}  \label{And}
\end{align}
This number is also the determinant of the periodic H\"uckel matrix,  $\det H_n(1,1)$.

A plane partition $\pi $ of $N$ of shape $(a,b,c)$ is an array $a\times b$ of numbers $c\ge \pi_{ij} \ge 0$ ($i=1,.. ,a$, $j=1,..,b$) with the property that all rows and columns are weakly decreasing, and the sum of all 
$\pi_{ij}$ is $N$.\\ 
These are two plane partitions of 38 in $(3,5,7)$:
\begin{align*}
\begin{array}{ccccc}
5 &5 & 3 & 2 & 2 \\
5 &5 & 2 & 1&0 \\
4 & 2 & 1 & 1&0
\end{array} \qquad 
\begin{array}{ccccc}
7 &6 & 6 & 6 & 1 \\
5 &3 & 2 & 1& 0\\
1 & 0 &0  & 0&0
\end{array}
\end{align*} 
A plane partition may be represented by stacking $\pi_{ij}$ cubes on each square cell $(i,j)$ in the rectangle $a\times b$. 
A stacking of cubes defines ascending staircases that correspond to sets of non-intersecting  
walks in a lattice. The beautiful Gessel-Viennot theorem \cite{Gessel85} counts them as determinants of certain
binomial matrices\footnote{see also the nice presentation by Aigner in \cite{Aigner}}.

It was later discovered \cite{Guy89} that the projection of the bounding box $a\times b\times c$ on a oblique plane normal
to (1,1,1) is an hexagon $H(a,b,c)$ with side-lengths $a, b, c, a, b, c$ (in cyclic order) and angles $\pi/3$ in a triangular lattice. Each stacking is one-to-one with a lozenge tiling of the hexagon. A lozenge is the union of two triangular unit cells of the foreground triangular lattice, and has angles $\pi/3$ and $2\pi/3$.\\ 
%
\begin{figure}
\begin{center}
\includegraphics*[width=5cm,clip=]{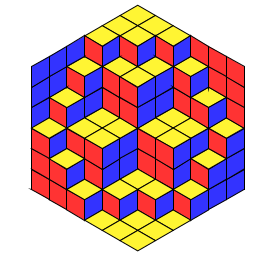}
\caption{\footnotesize A lozenge tile of $H(6,6,6)$ in Class 10. The cube stacking represents the partition with shape $6 \times 6 \times 6$ with rows:
6,6,6,5,4,3; 6,6,5,3,3,2; 6,5,5,3,3,1; 5,3,3,1,1,0; 4,3,3,1,1,0; 3,2,1,0,0,0 
(https://mathematica.stackexchange.com/questions/ 228422/lozenge-tilings)}
\end{center}
\end{figure}
Plane partitions were classified by Stanley and others into 10 symmetry classes and enumerated
(\cite{Krattenthaler15} is a nice review by Christian Krattenthaler).\\
Class 1 are the unrestricted partitions. They were enumerated by Percy MacMahon (1896):
the number of all plane partitions in $(a,b,c)$, i.e. lozenge tilings of $H(a,b,c)$, is
$$ \frac{H(a)H(b)H(c)H(a+b+c)}{H(a+c) H(a+c)H(b+c)} $$
where $H(k)=1! 2! ... (k-1)!$.\\
The other extremum is Class 10, with totally symmetric self-complementary plane partition in $(2n,2n,2n)$. 
Viewed as a stack
of cubes in the cubic box of side $2n$, a TSSC partition has the property of being equal to its complement (the void in the box,
see fig.3) and, viewed as a tiling, it is invariant under rotations by $\pi/3$ and reflections in the main diagonals of the hexagon. 
The enumeration was finally obtained by Andrews \cite{Andrews94}:
$$A(n) = \prod_{k=0}^{n-1} \frac{(3k+1)!}{(n+k)!} $$
$A(n)$ is also the number of $n\times n$ matrices $\{ 0,\pm 1 \}$ whose row-sums and column-sums are all 1,
and in every row and column the non-zero entries alternate in sign \cite{Zeilberger96}.

With the aid of the identity 
\begin{align}
1=\prod_{k=0}^n \frac{k!(n+k+1)!}{(2k)! (2k+1)!} 
\end{align}
the numbers \eqref{And} can be rewritten as follows:
\begin{align}
\det (Q_n+I_{n+1}) 
= A(n+1)\prod_{k=0}^n \frac{3k+2}{3k+1} 
\end{align}

Ciucu, Eisenk\"obl, Krattenthaler and Zare 
evaluated weighted sums over the tilings of the hexagon $H(a,b+m,c,a+m,b,c+m)$ with a central triangular
hole of sides $m,m,m$. With $a=b=c=n+1$ the counts are given by the determinant \eqref{dets} 
with weights $\omega =\pm 1,  e^{i2\pi/3},  e^{i\pi/3}$  (equations 3.3, 3.4, 3.5 in \cite{Ciucu99}). \\
With $m=0$ they simplify to eqs.\eqref{Ciucu3.3}, \eqref{Ciucu3.4} and \eqref{Ciucu3.5} given below. 
They are the weighted counts 
$\sum_\pi  \omega^{\# \pi}$ ($\# \pi $ is the number of cubes on the main diagonal) over the cyclically symmetric partitions $\pi $ with shape $(n+1)^3$
(Corollary 8 of \cite{Ciucu99}). In other words, the determinants give the weighted sums over the tilings of $H(n+1,n+1,n+1)$ that are invariant under rotations by $\pi/3$ around the center of the hexagon (Class 3). The different weights $\omega $ give:
\begin{align}
 \det (Q_{2n+1}-I_{2n+2} )= (-1)^{n+1}\left[ \prod_{k=0}^n \frac{k! (3k+1)!}{(2k)! (2k+1)!}\right]^4  
 = (-1)^{n+1} A(n+1)^4
 \label{Ciucu3.3}
\end{align}
Up to sign, it coincides with $\det H_{2n+1} (1,-1)$. For an odd size $\det (Q_{2n}-I_{2n+1})=0$. Some numbers are 
listed by Stembridge (1998), who first studied this case\footnote{page 6 in http://www.math.lsa.umich.edu/~jrs/papers/strange.pdf}.\\
Eqs. 3.4 and 3.5 with $m=0$ have simpler expressions given in Table 1 of \cite{Mitra09}:
\begin{align}
&\det (Q_n+e^{i\tfrac{2\pi}{3}} I_{n+1}) =  e^{i\frac{\pi}{3}(n+1)} A(n+1) \label{Ciucu3.4}\\
&\det (Q_n+e^{i\tfrac{\pi}{3}} I_{n+1}) =  e^{i\frac{\pi}{6}(n+1)} A_{HT}(n+1)^2 \times \begin{cases} 1 & \text{$n$ odd}\\
\sqrt 3 & \text{$n$ even} \end{cases}.   \label{Ciucu3.5}
\end{align}
The numbers $A_{HT}(n)$ enumerate the half-turn symmetric alternating sign $n\times n$ matrices\footnote{It is the
subclass of alternating sign $n\times n$ matrices such that $A_{ij}=A_{n+1-i, n+1-j}$. They are related to the ice model
of statistical mechanics.} \cite{Razumov05}:
\begin{align}
A_{HT}(2n) =  \prod_{k=0}^{n-1} \frac{(3k)!(3k+2)!}{[(n+k)!]^2} \qquad
A_{HT}(2n+1) = \frac{n! (3n)!}{[(2n)!]^2} A_{HT} (2n) 
\end{align}

The generalized Pascal matrices also occur in the works by Mitra and Nienhuis
\cite{MitraNienhuis04, Mitra09}. On an infinite cylindric square lattice of even circumference $L$, they consider 
coverings of closed paths that at each vertex make a left or right turn with equal probability. 
Two paths can have a vertex in common, but do not intersect.\\
The probability $P(L,m)$ that a point (not a vertex) of the cylinder is surrounded by $m$ loops is guessed as
$Q(L,m)/A_{HT}(L)^2$, where $Q(L,m)$ is expressed with the  
polynomial coefficients of the Pascal matrix. At the same time, $P(L,m)$ is estimated by Coulomb gas techniques.  
The following asymptotics is obtained:
\begin{align}
\det (Q_{L-1} + i I_L) &= i^{L/2}  L^{-5/48}  A_{HT}(L)^2 \times  \label{Mitra}\\
&\left [ 0.81099753 - \frac{0.028861}{L^{3/2}} +\frac{0.021012}{L^2}+ \frac{a_7}{L^{7/2}} + ...\right ] \nonumber
\end{align}
The large $n$ (i.e. $L+1$) expansion of \eqref{char} for any value $|\omega |=1$ is discussed by Mitra \cite{Mitra09}.\\

The H\"uckel matrices $H_n(\theta) = H_n(e^{-i\theta}, e^{i\theta})$ are Hermitian and their determinants can be written as
polynomials of $\cos\theta$. The known results for the Pascal matrix give:
\begin{gather*}
\det H_n(0) =  A(n+1) \prod_{k=0}^n \frac{3k+2}{3k+1} \\
 \det H_{2n+1} (\tfrac{\pi}{2}) = [A(n+1)]^4, \qquad \det H_{2n}(\tfrac{\pi}{2}) =0 \\
 \det H_n(\tfrac{\pi}{3}) = A(n+1)\\
 \det H_{2n+1}(\tfrac{\pi}{6}) = [A_{HT}(2n+2)]^2, \qquad  \det H_{2n}(\tfrac{\pi}{6}) = [A_{HT}(2n+1)]^2 \sqrt 3
\end{gather*}
The following table is a check. It shows the first values of $A(n)$, $A_{HT}(n)$ and $\det H_n(\theta)$, that coincide
with the values of the determinants listed in introduction. The value for $\theta=\pi/4$ is compared with the expansion \eqref{Mitra} 
by Mitra; the values are very close.
\begin{align*}
\begin{array}{c|c|c||c|c|c|c|c||c}
n & A & A_{HT} & \theta=0 & \pi/6 & \pi/3 &\pi/2 & \pi/4 & {\rm Mitra} \\ 
\hline
2 & 2 &2 &20&3^2\sqrt 3&7&0 & 8\sqrt 2 & \\
3 & 7 & 3 &132&10^2&42&2^4&70 &69.996 \\ 
4 & 42 &  10 &1452&25^2\sqrt 3&429 &0&526\sqrt 2 & \\
5 & 429 &25 &26741&140^2 &7436& 7^4&13167&13166.70 \\
6 & 7436 & 140 &826540&588^2\sqrt 3&218348&0&280772\sqrt 2 &\\
7 & 218348 & 588 &&&&& &\\   
\hline 
\end{array}
\end{align*}

\section{The second conjecture (honeycomb trapezia)}
If rows $0,1,...,k-1$ are removed from a honeycomb triangle, a trapezium results, with rows 
of lengths $2k+1, ... , 2n+1$.
The corresponding H\"uckel matrix $H_{k,n}({\bf x},{\bf y})$ is obtained by deleting the first $k^2$ rows and columns of $H_n({\bf x},{\bf y})$. Its size is
$(n+1)^2-k^2=(2k+1) + ... + (2n+1)$.

For example, the graph and the matrix $H_{6,7}$ (size $28$) are:\\
\begin{minipage}{.4 \textwidth}
\begin{center}
\begin{tikzpicture}
\draw[densely dotted, thick] (-3.5,-0.5) -- (-4,-0.5);
\filldraw[blue] (-3.5,-0.5) circle (2pt) node[anchor=north] {1};
\draw[gray, thick] (-3,0) -- (-3.5,-0.5);
\filldraw[red] (-3,0) circle (2pt) node[anchor=north] {};
\draw[gray, thick] (-3,0) -- (-2.5,-0.5);
\filldraw[blue] (-2.5,-0.5) circle (2pt) node[anchor=north] {};
\draw[gray, thick] (-2.5,-0.5) -- (-2,0);
\filldraw[red] (-2, 0) circle (2pt) node[anchor=north] {};
\draw[gray, thick] (-2,0) -- (-1.5,-0.5);
\filldraw[blue] (-1.5,-0.5) circle (2pt) node[anchor=north] {};
\draw[gray, thick] (-1.5,-0.5) -- (-1,0);
\filldraw[red] (-1,0) circle (2pt) node[anchor=north] {};
\draw[gray, thick] (-1,0) -- (-0.5,-0.5);
\filldraw[blue] (-0.5,-0.5) circle (2pt) node[anchor=north] {};
\draw[gray, thick] (-0.5,-0.5) -- (0,0);
\filldraw[red] (0,0) circle (2pt) node[anchor=north] {};
\draw[gray, thick] (0,0) -- (0.5,-0.5);
\filldraw[blue] (0.5,-0.5) circle (2pt) node[anchor=north] {};
\draw[gray, thick] (0.5,-0.5) -- (1,0);
\filldraw[red] (1,0) circle (2pt) node[anchor=north] {};
\draw[gray, thick] (1,0) -- (1.5,-0.5);
\filldraw[blue] (1.5,-0.5) circle (2pt) node[anchor=north] {};
\draw[gray, thick] (1.5,-0.5) -- (2,0);
\filldraw[red] (2,0) circle (2pt) node[anchor=north] {};
\draw[gray, thick] (2,0) -- (2.5,-0.5);
\filldraw[blue] (2.5,-0.5) circle (2pt) node[anchor=north] {};
\draw[gray, thick] (2.5,-0.5) -- (3,0);
\filldraw[red] (3,0) circle (2pt) node[anchor=north] {};
\draw[gray, thick] (3,0) -- (3.5,-0.5);
\filldraw[blue] (3.5,-0.5) circle (2pt) node[anchor=north] {15};
\draw[gray, thick] (-3,0) -- (-3,0.5);
\draw[gray, thick] (-2,0) -- (-2,0.5);
\draw[gray, thick] (-1,0) -- (-1,0.5);
\draw[gray, thick] (0,0) -- (0,0.5);
\draw[gray, thick] (1,0) -- (1,0.5);
\draw[gray, thick] (2,0) -- (2,0.5);
\draw[gray, thick] (3,0) -- (3,0.5);
\draw[densely dotted, thick] (-3,0.5) -- (-3.5,0.5);
\filldraw[blue] (-3,0.5) circle (2pt) node[anchor=south] {1};
\draw[gray, thick] (-3,0.5) -- (-2.5,1);
\filldraw[red] (-2.5,1) circle (2pt) node[anchor=south] {};
\draw[gray, thick] (-2.5,1) -- (-2,0.5);
\filldraw[blue] (-2, 0.5) circle (2pt) node[anchor=south] {};
\draw[gray, thick] (-2,0.5) -- (-1.5,1);
\filldraw[red] (-1.5,1) circle (2pt) node[anchor=south] {};
\draw[gray, thick] (-1.5,1) -- (-1,0.5);
\filldraw[blue] (-1,0.5) circle (2pt) node[anchor=south] {};
\draw[gray, thick] (-1,0.5) -- (-0.5,1);
\filldraw[red] (-0.5,1) circle (2pt) node[anchor=south] {};
\draw[gray, thick] (-0.5,1) -- (0,0.5);
\filldraw[blue] (0,0.5) circle (2pt) node[anchor=south] {};
\draw[gray, thick] (0,0.5) -- (0.5,1);
\filldraw[red] (0.5,1) circle (2pt) node[anchor=south] {};
\draw[gray, thick] (0.5,1) -- (1,0.5);
\filldraw[blue] (1,0.5) circle (2pt) node[anchor=south] {};
\draw[gray, thick] (1,0.5) -- (1.5,1);
\filldraw[red] (1.5,1) circle (2pt) node[anchor=south] {};
\draw[gray, thick] (1.5,1) -- (2,0.5);
\filldraw[blue] (2,0.5) circle (2pt) node[anchor=south] {};
\draw[gray, thick] (2,0.5) -- (2.5,1);
\filldraw[red] (2.5,1) circle (2pt) node[anchor=south] {};
\draw[gray, thick] (2.5,1) -- (3,0.5);
\filldraw[blue] (3,0.5) circle (2pt) node[anchor=south] {13};
\draw[densely dotted,thick] (3,0.5) -- (3.5,0.5);
\draw[densely dotted,thick] (3.5,-0.5) -- (4,-0.5);
\end{tikzpicture}
\end{center}

\end{minipage}
\begin{minipage}{.8 \textwidth}
\begin{align*}
\left [\begin{array}{cc} T_6 & R_7^T\\
R_7 & T_7 
\end{array}\right ]
\end{align*}
\end{minipage}

{\bf Conjecture 2:}
{\em The determinant of the H\"uckel matrix $H_{k,n}({\bf x}, {\bf y})$ of a honeycomb trapezium with rows 
with $2k+1$, $2k+3$, ... , $2n+1$ sites, size $(n+1)^2-k^2$, is equal to the determinant of the following matrix of size $n+1-k$}:
\begin{align*}
\det H_{k,n}({\bf x},{\bf y})= \left | \begin{array}{ccccc}
x_n+y_n & -\binom{n}{1} y_n & \binom{n}{2} y_n &  ... & \pm\binom{n}{n-k}y_n \\
\binom{n}{1}x_n & x_{n-1}+y_{n-1} & -\binom{n-1}{1}y_{n-1}  &  ... & \mp \binom{n-1}{n-k-1}y_{n-1}\\
\vdots &\vdots & \vdots &  ... & \vdots \\
\binom{n}{n-k} x_n &\binom{n-1}{n-k-1}x_{n-1}  & ... & ... & x_k+y_k\\
\end{array} \right |
\end{align*}

In these examples the determinants of H\"uckel matrices of size $28$, $51$ and $64$ are evaluated as matrices of size $2$, $3$ and $4$ ($S_n=x_n+y_n$):
\begin{align*} 
\det H_{6,7} = & \left | \begin{array}{cc} S_7 &  -7y_7\\ 7x_7& S_6\end{array} \right | =S_6 S_7 +7^2 x_7 y_7
\end{align*}
\begin{align*}
\det H_{7,9}=& \left | \begin{array}{ccc}
S_9 & -9 y_9 & 36 y_9  \\
9x_9 & S_8 & -8y_8 \\
36 x_9 &8x_8  &  S_7\\
\end{array} \right | = S_9 S_8 S_7 + 8^2x_8y_8 S_9+9^2 x_9y_9 S_7 + 36^2 x_9y_9 S_8
\end{align*}
\begin{align*}
\det H_{6,9}=& \left | \begin{array}{cccc}
S_9 & -9 y_9 & 36 y_9 & -84y_9  \\
9x_9 & S_8 & -8y_8 & 28 y_8\\
36 x_9 &8x_8  &  S_7 &-7y_7 \\
84 x_9 & 28 x_8 & 7 x_7 & S_6
\end{array} \right | 
= \, S_9[S_6 S_7 S_8 + 7^2 S_8 x_7 y_7 + 8^2 S_6 x_8 y_8 \\
&+ (28)^2 S_7 x_8 y_8] + x_9 y_9   [ 9^2 S_6 S_7    + 
 (36)^2 S_6 S_8]   + (63)^2 x_7 y_7 x_9y_9 \\
 & +(84)^2 (x_7 x_8 +y_7 y_8  +4 y_7 x_8 +4 x_7 y_8+16 x_8 y_8) x_9 y_9  
\end{align*}
It appears that all coefficients are perfect squares.

\section{The third conjecture (permanents)}
The equality of the permanent and the determinant of a matrix is a rare circumstance
\cite{Gibson71,McCuaig04}.
Numerical checks on H\"uckel matrices for triangles and trapezia of small sizes, 
support the conjecture:\\

{\bf Conjecture 3.}
{\em The permanent and the determinant of $H_n({\bf x},{\bf y})$ coincide. The same is true for $H_{k,n}({\bf x},{\bf y})$. }\\

Accordingly, $ \det H_n ={\sum}_\sigma  H_{1,\sigma_1}H_{2,\sigma_2} ... H_{N,\sigma_N} $ 
where $\sigma$ are even permutations of $1,...,N$, with $N=(n+1)^2$.
Each nonzero term contains $n(n+1)$ unit factors and $n+1$ factors $x_j$ or $y_k$ and
is visualized as $N$ arrows $1\to \sigma_1$ ... $N\to \sigma_N$ along the edges of the graph.
All vertices of the graph participate with one outgoing and one incoming arrow that
make closed oriented self-avoiding loops and dimers 
(two opposite arrows on the same edge). Such a configuration is a ``loop covering" $G$ 
of the graph. 

Theorem (Harari \cite{Harari}): {\em Let $G$ be a covering of the graph with dimers  and oriented
loops. If $\epsilon_G$ is the number of dimers and oriented loops of even length in $G$, and $p_G$ 
is the product of the values of the edges in $G$, then:} 
\begin{align}
\det H_n =\sum_G (-1)^{\epsilon_G} p_G
\end{align}
For honeycomb triangles and trapezia only even permutations contribute to the determinant, $\epsilon_G=1$.

\section{Some Proofs}
Some facts about determinants of H\"uckel matrices of graphene triangles and trapezia are now proven.  
\begin{prop} \label{PP1}
$\det H_n (x,y)$ is a homogeneous palindromic polynomial of degree $n+1$ in $x$ and $y$:
\begin{align*}
\det H_n(x,y) = (x^{n+1}+y^{n+1}) + c_1 (x^n y+xy^n) + c_2 (x^{n-1}y^2 + x^2 y^{n-1}) + \dots
\end{align*}
\begin{proof} Consider the diagonal matrix
$P_n={\rm diag}[(1)(1,\frac{1}{t} , 1)(1,\frac{1}{t},1, \frac{1}{t},1)  ...  ] $, of size
$1+3+...+(2n+1)=(n+1)^2$. It is $\det P_n=t^{-n(n+1)/2}$ and
\begin{align}
tP_n H_n(x,y)P_n= H_n (tx, ty)
\end{align}
In particular: $\det H_n (\det P_n)^2  t^{(n+1)^2}= \det H_n(tx,ty)$, i.e.  
\begin{align}
t^{n+1} \det H_n (x,y) = \det H_n (t x, t y)\label{result}
\end{align}
Therefore, the nonzero terms in the expansion of $\det H_n (x,y)$ are the products 
$H_{1i_1}...H_{N,i_N}$ that contain monomials $x^k y^{n+1-k}$.\\ 
This property and the symmetry $\det H_n(x,y)=\det H_n(y,x)$ imply the statement. The first coefficient is 1 because
$\det H_n(x,0)=x^{n+1}$.
\end{proof}
\end{prop}
As the result \eqref{result} only depends on the positions of the non-zero matrix elements and not on their values, every single 
product $H_{1,\sigma_1}...H_{N,\sigma_N}$ that does not contain a monomial of degree $n+1$ is identically zero. 

The graph of the triangle of graphene with $n+1$ rows has $b=2n+1$ blue dots at the end of the rows, $b'=\tfrac{1}{2}n(n-1)$ blue dots in bulk, $r=\tfrac{1}{2}n(n+1)$ red dots. The total number of vertices is $b+b'+r=(n+1)^2$.\\
Depending on the numbering of vertices, one has different matrix representations of the graph. 

\begin{prop} 
The coefficients of the expansion of $\det H_n ({\bf x},{\bf y})$ in monomials of degree $n+1$, are perfect squares of integers.
\begin{proof} 
Let us number the blue dots at the row ends from $1$ to $2n+1$, the remaining $b'$ blue dots from $2n+2$ to $b+b'$ and the 
red dots from $b+b'+1$ to $(n+1)^2$. The result is the block representation of the adjacency matrix:
\begin{align}
\tilde H_n= \left[ \begin{array}{c|cc} B & 0 & C \\ \hline 0 & O' & \rho \\  C^T & \rho^T & O_r \end{array} \right ] 
\end{align}
- $B$ has size $(2n+1)\times (2n+1)$ and elements in $\{0,x_i, y_i\}$,\\ 
- $C$ is the matrix $(2n+1)\times r$  of connections $(0,1)$ of row-ends with red dots,\\ 
- $\rho $ has elements $0,1$ that are the connections of inner blue with red dots,\\
- the $b'$ blue vertices are not inter-connected: this is the zero $b'\times b'$ matrix $O'$ ,\\ 
- the $r$ red vertices are not inter-connected: this is the zero $r\times r$ matrix $O_r$.\\ 
For example ($S_0=x_0+y_0$):\\
\begin{minipage}{.6 \textwidth}
\begin{align*}
\begin{array}{c||ccccc|c|ccc}
&0&1&2&3& 4&5&6&7&8\\
\hline\hline
0&S_0&&  &&    &  &1& & \\
1&&0 &y_1&&   & & 1&1 & \\
2&&  x_1&0&&   & & 1& &1 \\
3&&        & &0 &y_2      & & & 1& \\
4&&        & & x_2 &0 & & & &1 \\
\hline
5&& & & &       & 0 & & 1&1 \\
\hline
6&1&1 &1 & &    & &0 &0 &0 \\
7&&1 & &1 &   &1 &0 &0 & 0\\
8&& &1 & &   1&1 & 0& 0&0 \\
\end{array}
\end{align*}
\end{minipage}
\begin{minipage}{.5 \textwidth}
\begin{tikzpicture}
\draw[densely dotted,thick] (-2,-1) -- (-2.5,-1);
\draw[densely dotted,thick] (2,-1) -- (2.5,-1);
\draw[densely dotted,thick] (-1,0) -- (-1.5,0);
\draw[densely dotted,thick] (1,0) -- (1.5,0);
\draw[densely dotted,thick] (-0.5,1) -- (0.5,1);
\draw[gray, thick] (-2,-1) -- (2,-1);
\draw[gray, thick] (-1,-1) -- (1,-1);
\draw[gray, thick] (-1,-1) -- (-1,0);
\draw[gray, thick] (1,-1) -- (1,0);
\draw[gray, thick] (-1,0) -- (1,0);
\draw[gray, thick] (0,0) -- (0,1);
\filldraw[blue] (0,1) circle (2pt) node[anchor=south] {0};
\filldraw[blue] (1,0) circle (2pt) node[anchor=south] {2};
\filldraw[red] (0,0) circle (2pt) node[anchor=north] {6};
\filldraw[blue] (-1,0) circle (2pt) node[anchor=south] {1};
\filldraw[blue] (2,-1) circle (2pt) node[anchor=north] {4};
\filldraw[red] (1,-1) circle (2pt) node[anchor=north] {8};
\filldraw[blue] (0,-1) circle (2pt) node[anchor=north] {5};
\filldraw[red] (-1,-1) circle (2pt) node[anchor=north] {7};
\filldraw[blue] (-2,-1) circle (2pt) node[anchor=north] {3};
\end{tikzpicture}
\end{minipage}\\

The permutations of rows (1,2), (3,4) and (7,8) bring the
matrix to a block form that is diagonal in the parameters, associated to a graph that topologically differs
only in the end dots, that are now self-linked. The rest of the graph is unchanged by
the axial symmetry. As a consequence the new matrix $\hat H_n$ is symmetric: \\
\begin{minipage}{.7 \textwidth}
\begin{align*}
\begin{array}{c||ccccc|c|ccc}
&0&1&2&3& 4&5&6&7&8\\
\hline\hline
0&S_0&&  &&    &  &1& & \\
1&&  x_1&&&   & & 1& &1 \\
2&& &y_1&&   & & 1& 1& \\
3&&        & & x_2 & & & & &1 \\
4&&        & & &y_2      & & & 1& \\
\hline
5&& & & &       & 0 & & 1&1 \\
\hline
6&1&1 &1 & &    & &0 &0 &0 \\
7&& &1 & &   1&1 & 0& 0&0 \\
8&&1 & &1 &   &1 &0 &0 & 0\\
\end{array}
\end{align*}
\end{minipage}
\begin{minipage}{.5 \textwidth}
\begin{tikzpicture}
\draw[gray, thick] (-2,-1) -- (2,-1);
\draw[gray, thick] (-1,-1) -- (1,-1);
\draw[gray, thick] (-1,-1) -- (-1,0);
\draw[gray, thick] (1,-1) -- (1,0);
\draw[gray, thick] (-1,0) -- (1,0);
\draw[gray, thick] (0,0) -- (0,1);
\filldraw[black] (0,1) circle (2pt) node[anchor=south] {0};
\filldraw[black] (1,0) circle (2pt) node[anchor=south] {1};
\filldraw[red] (0,0) circle (2pt) node[anchor=north] {6};
\filldraw[black] (-1,0) circle (2pt) node[anchor=south] {2};
\filldraw[black] (2,-1) circle (2pt) node[anchor=north] {3};
\filldraw[red] (1,-1) circle (2pt) node[anchor=north] {8};
\filldraw[blue] (0,-1) circle (2pt) node[anchor=north] {5};
\filldraw[red] (-1,-1) circle (2pt) node[anchor=north] {7};
\filldraw[black] (-2,-1) circle (2pt) node[anchor=north] {4};
\end{tikzpicture}
\end{minipage}\\

According to Prop.\ref{PP1}, the determinant is a sum of homogeneous terms of degree $n+1$. An expansion by diagonal 
of the matrix gives:
$$ \det H_n ({\bf x},{\bf y}) = {\rm sign} \times {\sum}_{\boldsymbol \tau} S_0^{\tau_0} x_1^{\tau_1}...y_n^{\tau_{2n}} C_{\boldsymbol \tau}  $$
where $\tau_i =0,1 $ and $\sum_{i=0}^{2n} \tau_i =n+1$.  The coefficient $C_{\boldsymbol \tau} $ is a principal minor of 
$\hat H_n$ of size $(n+1)^2-n$, obtained by deleting the $n$ rows and $n$ columns that contain parameters with power 1, 
the other
parameters being put equal 0. The principal minor times sign is a squared determinant.\\
In the example, the sign is $(-1)^3$ and the coefficient (including sign) of $S_0x_1y_1$ is 
\begin{align*} 
-C_{1,1,1,0,0}=
-\left | \begin{array}{ccc|ccc}
0  & &  &  1& 0&1 \\
&  0&  & 1 &1 & 0\\
 &   &  0 & 0 &1 &1 \\
 \hline
1&1  &  0   & 0& 0&0 \\
 0 &  1  & 1  & 0& 0& 0\\
 1 & 0 & 1  & 0& 0& 0
\end{array}\right | =
&
 \left | \begin{array}{ccc}
    1& 0&1 \\
  1&1 &0 \\
  0 &1 &1 \\
\end{array}\right |^2 = 2^2
\end{align*}
\end{proof}
\end{prop}
\begin{remark} The signed principal minor has the meaning of matching number of the subgraph (Kasteleyn theorem). 
In the example the subgraph (in the left) admits only two coverings with dimers (i.e. matchings) that are shown.\\
\begin{center}
\begin{tikzpicture}
\draw[gray,thick] (-1,-1) -- (1,-1);
\draw[gray,thick] (-1,-1) -- (-1,0);
\draw[gray,thick] (1,-1) -- (1,0);
\draw[gray,thick] (-1,0) -- (1,0);
\filldraw[black] (1,0) circle (2pt) node[anchor=south] {1};
\filldraw[red] (0,0) circle (2pt) node[anchor=south] {6};
\filldraw[black] (-1,0) circle (2pt) node[anchor=south] {2};
\filldraw[red] (1,-1) circle (2pt) node[anchor=north] {8};
\filldraw[blue] (0,-1) circle (2pt) node[anchor=north] {5};
\filldraw[red] (-1,-1) circle (2pt) node[anchor=north] {7};
\end{tikzpicture}\qquad\qquad\quad
\begin{tikzpicture}
\draw[black, thick] (-1,-1) -- (0,-1);
\draw[black, thick] (1,-1) -- (1,0);
\draw[black, thick] (-1,0) -- (0,0);
\filldraw[black] (1,0) circle (2pt) node[anchor=south] {1};
\filldraw[red] (0,0) circle (2pt) node[anchor=south] {6};
\filldraw[black] (-1,0) circle (2pt) node[anchor=south] {2};
\filldraw[red] (1,-1) circle (2pt) node[anchor=north] {8};
\filldraw[blue] (0,-1) circle (2pt) node[anchor=north] {5};
\filldraw[red] (-1,-1) circle (2pt) node[anchor=north] {7};
\end{tikzpicture}\quad  \qquad\qquad 
\begin{tikzpicture}
\draw[black, thick] (0,-1) -- (1,-1);
\draw[black, thick] (-1,-1) -- (-1,0);
\draw[black, thick] (0,0) -- (1,0);
\filldraw[black] (1,0) circle (2pt) node[anchor=south] {1};
\filldraw[red] (0,0) circle (2pt) node[anchor=south] {6};
\filldraw[black] (-1,0) circle (2pt) node[anchor=south] {2};
\filldraw[red] (1,-1) circle (2pt) node[anchor=north] {8};
\filldraw[blue] (0,-1) circle (2pt) node[anchor=north] {5};
\filldraw[red] (-1,-1) circle (2pt) node[anchor=north] {7};
\end{tikzpicture}
\end{center}
\end{remark}

The same tricks used for honeycomb triangles apply to trapezia with simple modifications:
\begin{prop}\label{84}\quad\\
1) $\det H_{k,n} (x,y)$ is a homogeneous palindromic polynomial of degree $n+1-k$.\\
2) $\det H_{k,n} ({\bf x},{\bf y}) = (x_n+y_n) \det H_{k,n-1}({\bf x'},{\bf y'}) - x_n y_n \det H'_{k,n} ({\bf x'},{\bf y'})$,
where $H'_{k,n}$ is the matrix $H_{k,n}$ with the two rows and two columns of $x_n$,$y_n$ deleted.\\
3) If the parameters are distinct, the coefficients of the expansion of the determinant in monomials are perfect squares of
integers.
\end{prop}

\section{Analytic reduction for small matrices}
I now describe an analytic approach to the evaluation of the determinants. 
The key facts are the peculiar structure of the inverse of the diagonal blocks of $H_n$ and 
Schur's formula for determinants of block matrices.\\
Recall that $T_n$ has size $2n+1$.
 \begin{lem}\label{P3}
$R_n^T T_n^{-1} (x_n,y_n)R_n$ is a rank-1 matrix of size $2n-1$. For example:
\begin{align*}
& R_3^T T_3^{-1} (x_3,y_3)R_3= -\frac{x_3 y_3}{x_3+y_3} {\bf u}_3 {\bf u}_3^T,\quad {\bf u}_3^T = 
 \left[\begin{array}{cccccc}
+1 ,& 0, & -1, & 0,  & +1 \end{array}\right] \\
&R_4^T T_4^{-1} (x_n,y_n)R_4= \frac{x_4 y_4}{x_4+y_4} {\bf u}_4 {\bf u}_4^T, \quad {\bf u}_4^T = 
 \left[\begin{array}{ccccccc}
+1 ,& 0, & -1, & 0,  & +1, & 0,  & -1, 
\end{array}\right]
\end{align*}
\begin{proof}
The matrix $(x+y) T_n^{-1}(x,y)$ is the following:\\
$\bullet$ if $n=2k+1$ the matrix has size $4k+3$ and it is
\begin{align*}
\left [
\begin{array}{cccccccccccc}
(-1 & y & 1 & -y) & (-1 & y & 1 & -y )& \ldots  &  (-1 & y  & 1\\
x & -xy & y  & xy & -y &-xy & y  & 1 &    -y   &     & -xy & y \\
1 &x & -1 & y  & 1 & -y &1 & y  & 1  &  \ldots   & -y & -1 \\
-x &xy &x & -xy & y  & xy & -y &-xy & y  &     \ldots &  & -y  \\
-1 & -x & 1 & x & -1 & y & 1 & -y & 1 & \ldots  &  & 1\\
x & -xy & -x  &  &  &  &  &  & &           &   & \ldots \\
... &  &   &  &  &  &  &  &  &       &       &\ldots  \\
1 & x & -1  & -x &  &  &  &  &           & 1 & x  & -1
\end{array}\right ]
\end{align*}
$\bullet$ if $n=2k$ the matrix has size $4k+1$ and it is
\begin{align*}
\left[
\begin{array}{ccccccccccccc}
(1 & y & -1 & -y )& (1 & y & -1 & -y )& \ldots  &  & -1 & -y ) & 1\\
x & xy & y  & -xy & -y &xy & y  & -xy &     -y   & \ldots &    & -xy & -y \\
-1 &x & 1 & y  & -1 & -y &1 & y  & -1 &     -y   & \ldots &    & -1 \\
-x &-xy &x & xy & y  & -xy & -y &xy & y  & -xy    & \ldots &  & y  \\
1 & -x & -1 & x & 1 & y & -1 & -y & 1 & y  & -1 &   & 1\\
x & xy & -x  &  &  &  &  &  &  &       &  &    &\ldots \\
... &  &  &  &  &  &  &  &  &       &  &  & \ldots   \\
1 & -x &-1  & x  &  &  &  &  &  &       &    &  x & 1
\end{array}\right ] 
\end{align*}
(brackets are introduced to ease the pattern). If $xy=1$ the matrix is Toeplitz. \\
Multiplication of $T^{-1}$ with the rectangular matrices 
deletes the first and last rows and columns, and replaces even rows and columns with zeros.\\
For example: $R_3^T T_3^{-1} (x,y)R_3$ and $R_4^T T_4^{-1} (x,y)R_4$ are:
\begin{align*}
\frac{xy}{x+y} \left[\begin{array}{ccccc}
-1 & 0 & +1 & 0  & -1   \\
0 & 0 & 0  & 0 & 0        \\
+1 &0 & -1& 0  & +1  \\
0 &0 &0 & 0 & 0  \\
-1 & 0 & +1 &  0 & -1
\end{array}\right] ,\quad
\frac{xy}{x+y} \left[\begin{array}{ccccccc}
+1 & 0 & -1 & 0  & +1 & 0  & -1 \\
0  &  0 & 0  & 0 & 0    &0 &  0 \\
-1 &0 &+ 1& 0  & -1 & 0&  +1 \\
0 &0 &0 & 0 & 0  & 0 &   0\\
+1 & 0 & -1 &  0 & + 1&  0&  -1    \\
0 & 0 &  0 & 0& 0 &0  & 0   \\
  -1 &0 &+ 1& 0  & -1 & 0&  +1 
\end{array}\right]
\end{align*}
The matrices are then written in terms of vectors.
\end{proof}
\end{lem}
%

\begin{prop}\label{P4}
$\det H_n$ coincides with the determinant of the bordered matrix of size $n^2+1$:
\begin{align} 
\det H_n=\det \left [ \begin{array}{c|c} x_n+y_n & (-)^n x_n U^T\\
\hline
y_n U &H_{n-1} ({\bf x',y'})  \end{array} \right ]
\end{align}
where $U_n^T= [0,\ldots 0, {\bf u}_n]$ begins with $(n-1)^2$ zeros, 
followed by ${\bf u}_n$ that has $2n-1$ alternating components $1,0,-1,0,1,...$. 
\begin{proof}
Let us evaluate $\det H_n$ with Schur's formula :
\begin{align*}
\det \left[ \begin{array}{cc} A & B \\ C & D \end{array}\right ] = \det D \det (A-BD^{-1}C). 
\end{align*}
In this case, the full matrix is $H_n$, $D=T_n(x_n,y_n)$. With $\det T_n=x_n+y_n$, the formula and prop.\ref{P3} give:
\begin{align*}
\det H_n =(x_n+y_n) \det \left [H_{n-1} - (-1)^n\frac{x_n y_n}{x_n+y_n} U_n U_n^T \right ]
\end{align*}
If $a\neq 0$ is a number, Schur's formula states (Sherman-Morrison-Woodbury formula): 
\begin{align*}
\det \left[ \begin{array}{cc} a & {\bf b}^T \\ {\bf c} & M \end{array} \right ] = a \det (M-\frac{1}{a} {\bf c}{\bf b}^T) 
\end{align*}
With $a=x_n+y_n$, $M=H_{n-1}$, ${\bf b}^T= (-1)^n x_n U^T_n$ and ${\bf c} =y_n U_n$ the result follows.
\end{proof}
\end{prop}
%
These are the first bordered matrices:
\begin{gather*}
\det H_1=  \left | \begin{array}{c|c} 
 x_1+y_1 & -x_1 \\
 \hline
  y_1 & x_0+y_0 
  \end{array} \right | , \quad 
 \det H_2 = \left | \begin{array}{c|cccc} 
 x_2+y_2 & 0 & x_2 & 0 & -x_2 \\
 \hline
0 & x_0+y_0 & 0 & 1 & 0 \\
y_2& 0 & 0 & 1& y_1  \\ 
0& 1 & 1& 0&1  \\ 
-y_2& 0 & x_1& 1& 0 
\end{array} \right | 
\end{gather*}
\begin{align*}
\det H_3=\left | \begin{array}{c|ccccccccc}
x_3+y_3 &0 &0 &0 &0 & -x_3&0 &x_3 & 0& -x_3 \\
\hline
0&x_0+y_0 & 0 &1 &0 & 0& & & & \\
0&0& 0& 1& y_1& 0& 1& & & \\
0&1& 1& 0& 1&0 & & 0& & \\
0&0& x_1& 1& 0& 0& & & 1& \\
y_3&0&0 &0 &0 & 0& 1& 0& 0& y_2\\
0&& 1& & & 1& 0& 1& 0& 0\\
-y_3&& & 0& &0 & 1& 0& 1&0 \\
0&& & & 1& 0& 0& 1& 0& 1\\
y_3&& & & & x_2&0 &0 & 1& 0
\end{array}\right | 
\end{align*}

The bordering vectors and the off diagonal matrices $R_k$ have non-zero elements in different
rows and columns. This makes it possible to iterate the procedure with Schur's formula, that eliminates
a block but adds a border. In the end, the determinant of the matrix of size $(n+1)^2$ condensates 
to a bordered matrix of size $n+1$. This was done by hand for $n$ up to 5 to guess the binomial pattern
of the reduced matrix, and thus the connection with the Pascal matrix.

\section*{Aknowledgements}
I thank prof.~Stefano Evangelisti for addressing me to the study of determinants of H\"uckel matrices. This,
unexpectedly, turned me to the beautiful mathematics of Pascal matrices and combinatorics (that are in The Book \cite{Aigner}).  \\
I thank prof. Andrea Sportiello for his advices and insight for a possible proof.

\end{document}